# About some properties of simple trigonometric splines

Denysiuk V.P. Dr of Phys-Math. sciences, Professor, Kiev, Ukraine
National Aviation University
kvomden@nau.edu.ua

## Annotation

The class $Ts(r,f)$ the trigonometric interpolation splines depending on the parameter vectors, selected convergence factors and interpolation factors is considered. The main properties of simple interpolation trigonometric splines are given, which are also transferred to periodic simple interpolation polynomial splines. These results lead to the possibility of combining the theory of simple polynomial interpolation splines and the basics of the theory of simple trigonometric splines into a single theory - the theory of interpolation splines..

## Introduction

Approximation, respectively, the representation of a known or unknown function through a set of some special functions can be considered as a central topic of analysis; such special functions are well defined, easy to calculate, and have certain analytical properties [1]. Algebraic and trigonometric polynomials, exponential functions [2], polynomial [3] and trigonometric [4] splines, etc. often act as special functions.

It is known [5,6] that the best apparatus for approximating the functions of a class. $W_v^r$ there are simple polynomial splines where $W_v^r$ - a class of periodic functions that have an absolutely continuous derivative of order $r-1$ ( $r=1,2,...$ ), and a derivative of the order $r$ is a function of limited variation. The theory of such splines is well developed (look [3], [7], [8], [9]) etc. The main disadvantage of polynomial splines, in our opinion, is that they have a lumpy structure; this leads to the fact that in practice mainly splines of the third degree are used, which are sewn together from pieces of algebraic polynomials of the third degree. However, this structure of splines significantly limits their application in many problems of computational mathematics.

By Schoenberg [10], has been proposed that splines that are crosslinked from trigonometric polynomials of a certain order and have the same disadvantage as polynomial splines, a piecewise structure.

In [11] a class of generalized trigonometric functions was introduced, on the basis of which another principle for constructing trigonometric splines was proposed in [4] - their representation by uniformly convergent trigonometric series (Fourier series), the coefficients of which have a certain order of decline. The undeniable advantage of such splines is that they are given as a single expression over the entire interval of the function. Subsequently, such splines were considered in [12], [13], [14].

In the process of studying trigonometric splines, it was found that they have a number of properties similar to the properties of simple polynomial splines. In our opinion, this allows us to consider classes of trigonometric and classes of polynomials as a single class - the class of interpolation splines. This combination looks very promising because it allows you to combine research methods common to both classes.

## The purpose of the work

Investigation and comparison of properties of interpolation periodic simple polynomial and trigonometric splines.

## The main part

We introduce some concepts and notations. Following [Alberg], we denote by $K^n$ function class $f(t)$ , $t \in [0, 2\pi]$, which has a completely continuous derivative ($n-1$)-st order and $n$-th derived from space $L_{2\pi}^2$. Also class $K_p^n$ class $2\pi$ - periodic functions having a completely continuous derivative ($n-1$)-st order and $n$-th derived from space $L_{2\pi}^2$. Note that the class is introduced in this way $K_p^n$ coincides with the previously entered class $W_v^r$ while $n = r$.

It is clear that using special methods of periodic continuation, given in [15], the functions of the class $K^n$ can be translated into class functions $K_p^n$; therefore, in the future, without losing generality, we will confine ourselves to considering the functions of the class $K_p^n$. Note that often for classes $K^n$ and $K_p^n$ also use other designations. For functions $f(t)$ of class $K_p^n$ you can enter the half-norm [Alberg] as follows

$$\|f\| = \left\{ \int_a^b |f^{(n)}(t)|^2 \, dx \right\}^{\frac{1}{2}}.$$

In the theory of trigonometric series, another norm is often considered, which we will denote $P(f,k)$ and which is given by the expression

$$\|f\|_k = \frac{1}{\pi} \int_0^{2\pi} \left[ f^{(k)}(t) \right]^2 dx,$$

It is clear that the norm $\|f\|_k$ in a sense it is proportional to the half-norm $\|f\|$ and also characterizes the value $f^{(k)}(t)$. In applications, the value of $\|f\|_k$ often interpreted as the average power of the function $f^{(k)}(t)$ on the period; therefore further norm $\|f\|_k$ we will denote as $P(f,k)$.

Since the trigonometric splines introduced in [4], [12], [13], [14] are generated by interpolation trigonometric polynomials, it is expedient to begin consideration of some properties of trigonometric polynomials constructed in this way..

Let's $\Delta_N^{(I)} = \{t_i\}_{i=1}^N$ - uniform grids set on $[0, 2\pi)$, where $t_i^{(0)} = 2\pi \frac{i-1}{N}$, $t_i^{(1)} = \pi \frac{2i-1}{N}$, $N = 2n+1$, $I$ - indicator, $I = 0,1$. Let also on $[0, 2\pi)$ given some continuous function $f(t)$ and let $f_i^{(I)} = f\left(t_i^{(I)}\right)$ - the value of this function in the grid nodes $\Delta_N^{(I)}$. Then the interpolation trigonometric polynomial $T_n^{(I)}(f,t)$, interpolating function $f(t)$ in grid nodes $\Delta_N^{(I)}$, can be written as

$$T_n^{(I)}(t) = \frac{a_0}{2} + \sum_{k=1}^n a_k^{(I)} \cos kt + b_k^{(I)} \sin kt,$$

where the coefficients $a_0, a_k^{(I)}, b_k^{(I)}$ ($k = 1, 2, ..., n$) calculated by formulas [16]

$$a_k^{(I)} = \frac{2}{N} \sum_{j=1}^N f_j^{(I)} \cos kt_j^{(I)}, \quad b_k^{(I)} = \frac{2}{N} \sum_{j=1}^N f_j^{(I)} \sin kt_j^{(I)} \quad . \quad (3)$$
$$k = 0, 1, ..., n; \quad\quad\quad\quad k = 1, 2, ..., n.$$

Note that to calculate the coefficients $a_0, a_k^{(I)}, b_k^{(I)}$ ($k = 1, 2, ..., n$) are often used algorithms FFT (**fast Fourier transform**).

Because the value of the coefficient $a_0$ does not depend on the value of the indicator $I$, we omit this indicator in the expression for $T_n^{(I)}(f,x)$.

Interpolation trigonometric polynomial $T_n^{(I)}(f,t)$ can be submitted in another form [17]. We consider fundamental trigonometric functions $tm_k^{(I)}(t)$, which satisfy on the grid $\Delta_N^{(I)}$ the conditions[18]

$$tm_k^{(I)}\left(t_j^{(I)}\right) = \begin{cases} 1, & k = j \\ 0, & k \neq j \end{cases}, \quad (k, j = 1, ..., N). \quad (\Phi)$$

Such polynomials can be represented as [19]

$$tm_k^{(I)}(t) = \frac{1}{N} \left[ 1 + 2 \sum_{j=1}^n \cos j\left(t - t_k^{(I)}\right) \right].$$

Using a system of fundamental trigonometric polynomials $t_k(t)$, $k = 1, 2, \cdots, N$, the interpolation trigonometric polynomial can be written as

$$T_n^{(I)}(t) = \sum_{k=0}^{N} f_k^{(I)} \, tm_k^{(I)}(t).$$

Note that the fundamental trigonometric functions $t_k^{(I)}(t)$ satisfy the condition

$$\sum_{k=1}^{N} tm_k^{(I)}(t) = 1.$$

In addition, these functions are orthogonal in the sense that

$$\int_0^{2\pi} tm_k^{(I)}(t) \, tm_j^{(I)}(t) \, dx = \begin{cases} \pi, & k = j; \\ 0, & k \neq j. \end{cases}$$

Taking into account given the above, Parseval's theorem for an interpolation trigonometric polynomial $T_n^{*(I)}(f,t)$ can be submitted in the following forms.

$$\frac{1}{\pi} \int_0^{2\pi} \left[ T_n^{(I)}(t) \right]^2 dx = \frac{(a_0)^2}{2} + \sum_{k=1}^{n} (a_k^{(I)})^2 + (b_k^{(I)})^2 =$$

$$= \frac{2}{N} \sum_{j=1}^{N} \left( T_n^{(I)}(t_j^{(I)}) \right)^2 = \frac{2}{N} \sum_{j=1}^{N} \left( T_n^{(I)}(t_j^{(1-I)}) \right)^2 = \frac{2}{N} \sum_{j=1}^{N} \left( f_j^{(I)} \right)^2$$

Thus, we can determine the average power of the function or through the value of the function $f^{(I)}(t)$ in grid nodes $\Delta_N^{(I)}$, or due to the value of the trigonometric interpolation polynomial in the nodes of the same grids. Generally speaking, due to the periodicity of the trigonometric interpolation polynomial, to calculate the average power, you can use its value in the nodes of arbitrary uniform grids on $[0, 2\pi)$, containing $N$ grids.

Let us now turn to the consideration of trigonometric splines. Given [14], we present trigonometric interpolation splines on the grids in General as follows

$$St^{(I_1,I_2)}(\Gamma,H,\nu,r,t) = \frac{a_0^{(I_2)}}{2} + \sum_{k=1}^{\frac{N-1}{2}} \left[ a_\kappa^{(I_2)} \frac{C_k^{(I_1)}(\Gamma,\nu,r,t)}{hc_k^{(I_1,I_2)}(\Gamma,r,k)} + b_\kappa^{(I_2)} \frac{S_k^{(I_1)}(H,\nu,r,t)}{hs_k^{(I_1,I_2)}(H,r,k)} \right],$$

where

$$C_k^{(I_1)}(\Gamma,\nu,r,t) =$$
$$= \gamma_1 \nu_k(r) \cos kt + \sum_{m=1}^{\infty} (-1)^{mI_1} \left[ \gamma_2 \nu_{mN+k}(r) \cos(mN+k)t + \gamma_3 \nu_{mN-k}(r) \cos(mN-k)t \right];$$

$$S_k^{(I_1)}(H,\nu,r,t) = \tag{3}$$
$$= \eta_1 \nu_k(r) \sin kt + \sum_{m=1}^{\infty} (-1)^{mI_1} \left[ \eta_2 \nu_{mN+k}(r) \sin(mN+k)t - \eta_3 \nu_{mN-k}(r) \sin(mN-k)t \right],$$

with convergence factors $\nu_k(r)$, having a descending order $O(k^{-(1+r)})$, and interpolation factors

$$hc_k^{(I_1,I_2)}(\Gamma,\nu,r) = \gamma_1 \nu_k(r) + \sum_{m=1}^{\infty} (-1)^{m(I_1-I_2)} \left[ \gamma_2 \nu_{mN+k}(r) + \gamma_3 \nu_{mN-k}(r) \right],$$

$$hs_k^{(I_1,I_2)}(H,\nu,r) = \eta_1 \nu_k(r) + \sum_{m=1}^{\infty} (-1)^{m(I_1-I_2)} \left[ \eta_2 \nu_{mN+k}(r) + \eta_3 \nu_{mN-k}(r) \right]. \tag{4}$$

Where indicator $I_1$ ($I_1 = 0,1$) determines the stitching grid, indicator $I_2$ ($I_2 = 0,1$) determines the interpolation grid, $a_0^{(I_2)}, a_\kappa^{(I_2)}, b_\kappa^{(I_2)}$ - coefficients of the interpolation trigonometric polynomial on the grid $\Delta_N^{(I_2)}$, $r$, ($r = 1, 2, ...$) - parameter that determines the order of the spline, and $\Gamma = \{\gamma_1, \gamma_2, \gamma_3\}$ and $H = \{\eta_1, \eta_2, \eta_3\}$ - parameter vectors, and parameters $\gamma_k$ and $\eta_k$, ($k = 1,2,3$) take arbitrary real values and at least three of them do not rotate simultaneously by 0. Note that to reduce the notation of splines and the functions through which they are built, we omit the dependence on the number $N$ nodes of interpolation grids $\Delta_N^{(I)}$.

It is clear that the spline $St^{(1,1)}(\Gamma,H,\nu,r,t)$ is a shifted half-step grid $\Delta_N^{(I)}$ spline $St^{(0,0)}(\Gamma,H,\nu,r,t)$; so is the spline $St^{(0,1)}(\Gamma,H,\nu,r,t)$ is a spline shifted half a step of the same grid $St^{(1,0)}(\Gamma,H,\nu,r,t)$. However, in the following we will consider all four such splines.

Note that when $\Gamma = \{1,0,0\}$ and $H = \{1,0,0\}$ we have a regular trigonometric interpolation polynomial $T_n^{(I)}(f,x)$.

In the accepted designations it is clear that $St^{(I_1,I_2)}(\Gamma,H,\nu,r,t) \in K_p^r$.

If $\Gamma = \{1,1,1\}$ and $H = \{1,1,1\}$, then we have simple trigonometric splines, in the notation of which we will omit the dependence on the vectors $\Gamma$ and $H$. Simple trigonometric splines $St^{(0,0)}(\nu,r,t)$ with convergence factors

$$\nu 1_k(r) = \left[\sin c(\frac{\pi k}{N})\right]^{1+r}, \quad \nu 2_k(r) = \left[\left|\sin c(\frac{\pi k}{N})\right|\right]^{1+r}, \quad \nu 3_k(r) = \left[\frac{1}{k}\right]^{1+r}, \quad \nu 4_k(r) = sign\left[\sin(\frac{\pi k}{N})\right]\left[\frac{1}{k}\right]^{1+r},$$

and with $r = 2k-1$ ($k=1,2,...$) coincide with simple periodic polynomial splines of degree $r$. Generally speaking, with convergence factors $\nu 1_k(r)$ and $\nu 4_k(r)$ and when $r = 1, 2, ...$ all simple trigonometric splines have polynomial analogues, but most of these analogues are unknown.

As we saw above, the representation of trigonometric polynomials through fundamental polynomials has a number of advantages; therefore, we present trigonometric splines through fundamental splines. It is clear that fundamental trigonometric splines exist only if $\Gamma = H$; while fundamental splines $st_k^{(I_1,I_2)}(\Gamma,\nu,r,j,t)$ have the form

$$st_k^{(I1,I2)}(\Gamma,\nu,r,t) = \frac{1}{N}\left[1 + 2\sum_{k=1}^{\frac{N-1}{2}} \frac{c_k^{(I1)}(\Gamma,\nu,r,t)}{h_k^{(I1,I2)}(\Gamma,\nu,r)}\right],$$

where

$$c_k^{(I_1)}(\Gamma,\nu,r,t) = \gamma_1 \nu_k(r)\cos k(t - x_k^{(I_2)}) +$$

$$+ \sum_{m=1}^{\infty}(-1)^{mI_1}\left[\gamma_2 \nu_{mN+k}(r)\cos(mN+k)(t-x_k^{(I_2)}) + \gamma_3 \nu_{mN-k}(r)\cos(mN-k)(t-x_k^{(I_2)})\right];$$

$$hc_k^{(I_1,I_2)}(\Gamma,\nu,r) = \gamma_1 \nu_k(r) + \sum_{m=1}^{\infty}(-1)^{m(I_1-I_2)}\left[\gamma_2 \nu_{mN+k}(r) + \gamma_3 \nu_{mN-k}(r)\right],$$

Using fundamental splines $st_k^{(I1,I2)}(\Gamma,\nu,r,t)$, $k = 1,2,\cdots,N$, trigonometric spline $St^{(I_1,I_2)}(\Gamma,\nu,r,t)$ can be submitted as

$$St^{(I_1,I_2)}(\Gamma,\nu,r,t) = \sum_{k=1}^{N} f_k st_k^{(I1,I2)}(\Gamma,\nu,r,t).$$

Since in the theory of polynomial splines the derivatives of these splines are often considered, we introduce the derivatives of trigonometric splines. Such derivatives are of the order $q$, ($q \leq r-1$) splines $St^{(I_1,I_2)}(\Gamma,\nu,r,t)$ and fundamental splines $st_k^{(I1,I2)}(\Gamma,\nu,r,t)$ ($k = 1,2,...,N$) we will denote accordingly $St^{(I_1,I_2)}(\Gamma,\nu,r,q,t)$ and $st_k^{(I1,I2)}(\Gamma,\nu,r,q,t)$. So we have

$$\left[St^{(I_1,I_2)}(\Gamma,\nu,r,t)\right]^{(q)} = St^{(I_1,I_2)}(\Gamma,\nu,r,q,t);$$

$$\left[\sum_{k=1}^{N} f_k st_k^{(I1,I2)}(\Gamma,\nu,r,t)\right]^q = \sum_{k=1}^{N} f_k st_k^{(I1,I2)}(\Gamma,\nu,r,q,t).$$

In addition, further simple periodic polynomial splines crosslinked from degree polynomials $2m-1$, we will note as $Sp_{2m-1}(t)$; it is clear that $Sp_{2m-1}(t) \in K_p^m$.

Here are some properties of trigonometric splines.

1. 1. In the vast majority of cases, functions approaching simple polynomial splines are not periodic. However, as we have already said, using special methods of periodic continuation [15], they can be extended periodically to the entire numerical axis, and hence to approximate them by periodic polynomial or trigonometric splines of certain orders.
2. Simple spline $St^{(1,1)}(\nu,r,t)$ and its derivatives of order $q$ ($q=1,2,...,r-1$) are shifted by half a step grid $\Delta_N^{(I)}$ a simple spline $St^{(0,0)}(\nu,r,t)$ and its derivatives of the same order; just as simple a spline

$St^{(0,1)}(v,r,t)$ and its derivatives of the order are shifted by a half-step of the same grid a simple spline $St^{(1,0)}(v,r,t)$ and its derivatives of the same order.

3. Derivatives of order $r$ ($r=1,2,...$) simple trigonometric splines $St^{(I1,I2)}(v,r,t)$ coincide to piecewise - constant functions (it is necessary to consider the fact that these derivatives do not have properties of uniform convergence), and derivatives of the order $r-1$ of the same splines coincide evenly to periodic simple polynomial splines of the 1st order - broken, sewn together in grid nodes $\Delta_N^{(I1)}$. From this, in particular, follows the fact that simple trigonometric splines of order coincide uniformly with periodic functions, which are sewn together from pieces of algebraic polynomials with ensuring the continuity of derivatives to $r-1$-st order.

   Moreover, there is reason to believe that polynomial analogues of trigonometric splines with arbitrary vectors $\Gamma$ and $H$ also consist of algebraic polynomials that are sewn together in grid nodes $\Delta_N^{(I1)}$.

4. Polynomial analogues are known for trigonometric splines $St^{(0,0)}(v,2k-1,t)$ ($k=1,2,...$) - these are periodic simple polynomial splines of odd degree - broken, cubic, etc. [3]. Polynomial analogues are also known for the spline $St^{(0,1)}(v,2k,t)$; these are periodic simple polynomial splines of odd degree - broken, cubic are periodic simple polynomial splines of even degree [9]. For simple trigonometric splines $St^{(1,0)}(v,2k-1,t)$, $St^{(0,1)}(v,2k-1,t)$, $St^{(1,1)}(v,2k-1,t)$ odd order and $St^{(0,0)}(v,2k,t)$, $St^{(0,1)}(v,2k,t)$ $St^{(1,1)}(v,2k,t)$ in pairs, their polynomial analogues, although they exist, are still unknown.

5. From the identical equality of periodic simple polynomial splines $Sp_{2k-1}(t)$ and $Sp_{2k}(t)$ and simple trigonometric splines $St^{(0,0)}(v,2k-1,t)$ and $St^{(0,1)}(v,2k,t)$ it follows that all the results of the theory of approximations obtained for periodic simple polynomial splines can be transferred to simple trigonometric splines $St^{(0,0)}(v,2k-1,t)$ and $St^{(0,1)}(v,2k,t)$. Moreover, given property 2, these results can be transferred to splines $St^{(1,1)}(v,k+1,t)$ and $St^{(1,0)}(v,k+1,t)$

6. To construct periodic simple polynomial splines $Sp_{2k-1}(t)$ of order $2k-1$ it is necessary to find their moments, which are the values of their derivatives of order $2k-2$, ($k=2,3,...$) in grid nodes $\Delta_N^{(0)}$; to find these points we have to solve systems of algebraic equations of order $N$.

7. When constructing trigonometric splines, there is no need to find their moments. Moreover, on the contrary, the moments of periodic simple polynomial splines can be found as values of derivatives of order $2k-2$ trigonometric splines of order $2k-1$ in grid nodes $\Delta_N^{(0)}$; it is clear that algorithms FFT can be widely used.

8. In [4] it was shown that

$$\lim_{r \to \infty} St^{(I_1,I_2)}(\Gamma,H,v,r,t) = T_n^{(I2)}(t).$$

9. From properties (3), (4) and (7) it follows that

$$\lim_{r \to \infty} Sp_r(t) = T_n^{(0)}(t).$$

10. Since trigonometric splines are actually represented by their Fourier series, they have a Parseval theorem, which has the form

$$\frac{1}{\pi} \int_0^{2\pi} \left[ St^{(I_1,I_2)}(\Gamma,H,v,r,t) \right]^2 = \frac{(a_0)^2}{2} + pc^{(I1,I2)} + ps^{(I1,I2)},$$

where

$$pc^{(I1,I2)}(\Gamma,v,r) = \sum_{k=1}^{\frac{N-1}{2}} \left\{ \left( \frac{a_k^{(I2)}}{hc_k^{(I1,I2)}(\Gamma,v,r)} \right)^2 \left[ \left(\gamma_1 V_k(r)^2\right) + \sum_{m=1}^{\infty} \left[ \left(\gamma_2 V_{mN+k}(r)\right)^2 + \left(\gamma_3 V_{mN-k}(r)\right)^2 \right] \right] \right\};$$

$$ps^{(I1,I2)}(\mathrm{H},v,\mathrm{r}) = \sum_{k=1}^{\frac{N-1}{2}} \left\{ \left( \frac{b_k^{(I2)}}{hs_k^{(I1,I2)}(\Gamma,v,r)} \right)^2 \left[ \left(\eta_1 v_k(\mathrm{r})^2\right) + \sum_{m=1}^{\infty} \left[ \left(\eta_2 v_{mN+k}(\mathrm{r})\right)^2 + \left(\eta_3 v_{mN-k}(\mathrm{r})\right)^2 \right] \right] \right\}.$$

11. From property (7) it follows that

$$\lim_{r \to \infty} P\left(St^{(I1,I2)}(\Gamma,\mathrm{H},v,r),q\right) = P\left(T_n^{(I2)},q\right), \qquad (q=0,1,...,r-1)$$

12. From property (7) it follows that

$$\lim_{r \to \infty} P\left(Sp_{2r-1},q\right) = P\left(T_n^{(0)},q\right), \qquad (q=0,1,...,r-1).$$

13. Properties (7), (8), (10), (11) also occur when $N \to \infty$.

14. In case when $\Gamma = \mathrm{H}$, equalities under take place
$$P\left(St^{(0,0)}(\Gamma,v,r),q\right) = P\left(St^{(1,1)}(\Gamma,v,r),q\right);$$

$$P\left(St^{(0,1)}(\Gamma,v,r),q\right) = P\left(St^{(1,0)}(\Gamma,v,r),q\right), \qquad (q=0,1,...,r-1).$$

Given these equations, it is clear that in the study of quantities $P(f,q)$ can be limited to the consideration of trigonometric splines $St^{(0,0)}(\Gamma,v,r)$, $St^{(1,0)}(\Gamma,v,r)$ and their derivative of order $q$, ( $q \le r-1$).

15. It is known [3] that of all functions $f(t) \in K_p^m$, taking the specified values on the grid $\Delta_N^{(0)}$, minimum value $P(f,m)$ is achieved on simple polynomial interpolation splines $Sp_{2m-1}(t)$; this property is called the property of minimum curvature. However, it does not follow that the spline itself or its derivatives are of the order $q$ ($q=1,2,...,k-1$) have the same property. Indeed, it is easy to choose parameter vectors $\Gamma$ and $\mathrm{H}$, under which the condition will be fulfilled
$$P\left[St^{(I_1,I_2)}(\Gamma,\mathrm{H},v,2k-1),q\right] < P\left[Sp_{2k-1},q\right], \qquad (k=2,3,...).$$

16. Fundamental trigonometric splines exist only when $\Gamma = \mathrm{H}$.

17. It is easy to verify that the fundamental trigonometric splines satisfy the condition
$\sum_{k=1}^{N} st_k^{(I1,I2)}(\Gamma,v,r,0,t) = 1$ for all $\Gamma$.

18. It is easy to verify that the derivatives are in order $q \le r-1$ simple fundamental trigonometric splines satisfy the condition
$$\sum_{k=1}^{N} st_k^{(I1,I2)}(\Gamma,v,r,q,t) = 0 \text{ for all } \Gamma.$$

Properties (14), (15) of fundamental trigonometric splines coincide with similar properties of fundamental trigonometric functions.

19. Fundamental trigonometric splines, in contrast to fundamental trigonometric functions, do not have the properties of orthogonality. From this, in particular, it follows that in the study of the average power of trigonometric splines in the period, their representation through the fundamental splines is impractical.

# Conclusion

The main properties of simple interpolation trigonometric splines are given, which are easily transferred to simple interpolation polynomial splines. This possibility of transfer enriches both the theory of simple interpolation trigonometric splines and the theory of simple interpolation polynomial splines; moreover, it leads to the possibility of combining both theories into a single one - the theory of simple interpolation splines. Of course, it should be borne in mind that the set of simple trigonometric splines is a subset of the set of trigonometric splines.

Of course, the considered theory of simple interpolation splines needs further research.

**List of references**